\documentclass[12pt,a4paper,reqno]{amsart}
\usepackage{amsmath}
\usepackage{amsfonts}
\usepackage{amssymb}
\numberwithin{equation}{section}

\theoremstyle{plain}

\newtheorem{theorem}[subsection]{Theorem}

\newtheorem{lemma}[subsection]{Lemma}
\newtheorem{corollary}[subsection]{Corollary}
\newtheorem{conjecture}[subsection]{Conjecture}

\newtheorem*{quote-thm-1}{Proposition 8.1 of \cite{green-tao-longprimeaps}}
\newtheorem*{mainthm-restate}{Theorem \ref{mainthm-dual}}
\newtheorem*{convex-lemma-restate}{Lemma \ref{convex-props}}
\newtheorem*{decomp-prop-restate}{Proposition \ref{decomp-prop}}
\newtheorem*{main-theorem}{Main Theorem}
\newtheorem*{gvn-restate}{Proposition \ref{gvn}}
\newtheorem*{gvn-restate-again}{Proposition $\mbox{\ref{gvn}}^{\prime}$}
\newtheorem*{gvn-restate-yetagain}{Proposition $\mbox{\ref{gvn}}^{\prime \prime}$}
\newtheorem*{pseudodom-repeat}{Proposition \ref{pseudodom}}

\theoremstyle{definition}

\newtheorem{definition}[subsection]{Definition}

\theoremstyle{remark}

\renewcommand{\leq}{\leqslant}
\renewcommand{\geq}{\geqslant}

\newsavebox{\proofbox}
\savebox{\proofbox}{\begin{picture}(7,7)%
  \put(0,0){\framebox(7,7){}}\end{picture}}

\newcommand\GI{\operatorname{GI}}
\newcommand\MN{\operatorname{MN}}

\def\Z{\mathbb{Z}}
\def\R{\mathbb{R}}

\def\Q{\mathbb{Q}}

\newcommand\mobname{M\"obius and nilsequences\ }

\begin{document}

\title[Sublattices of finite index]{Sublattices of finite index}

\author{Chunlei Liu}
\address{School of Mathematics, Beijing Normal University, Beijing 100875}
\email{clliu@bnu.edu.cn} 
\thanks{This work is supported by Project 10671015 of the Natural Science Foundation of China.}

\begin{abstract}Assuming the Gowers Inverse conjecture and the
M\"{o}bius conjecture for the finite parameter $s$, Green-Tao
verified Dickson's conjecture for lattices which are ranges of
linear maps of complexity at most $s$. In this paper, we reformulate
Green-Tao's theorem on Dickson's conjecture, and prove that, if $L$
is the range of a linear map of complexity $s$, and $L_1$ is a
sublattice of $L$ of finite index, then $L_1$ is the range of a
linear map of complexity $s$.
\end{abstract}
\setcounter{tocdepth}{1}

\maketitle

\section{Introduction}\label{intro-sec}

A \emph{lattice} in an Euclidean vector space is an additive
discrete subgroup other than $\{0\}$. For example, $q\Z$ ($0\neq
q\in\Z$) is a sublattice of $\Z$, the set
$$\Z(1,1)=\{(n,n) \mid n\in\Z\}$$ is a sublattice of $\Z^2$, the set
$$\Z(1,0,-1)+\Z(0,1,-1)=\{(n_1,n_2,-n_1-n_2) \mid n_1,n_2\in\Z^2\}$$
is a sublattice of $\Z^3$, and
$$\Z(1,1,\cdots,1)+\Z(0,1,\cdots,t-1)=\{(n,n+r,\cdots,n+(t-1)r) \mid n,r\in\Z\},$$
which is the set of vectors $(n_1,\cdots,n_t)$ in $\Z^t$ such that
$n_1,\cdots,n_t$ form an arithmetic progression, is a sublattice of
$\Z^t$.

An \emph{affine sublattice} of $\Z^t$ is a translation of a
sublattice of $\Z^t$. For example, $q\Z+a$ ($q,a\in\Z, q\neq0$),
which is an arithmetic progression, is an affine sublattice of $\Z$,
the set
$$\Z(1,1)+(0,2)=\{(n,n+2) \mid n\in\Z\}$$ is an affine sublattice of $\Z^2$,
and the set
$$\Z(1,0,-1)+\Z(0,1,-1)+(0,0,N)=\{(n_1,n_2,N-n_1-n_2) \mid n_1,n_2\in\Z^2\}$$
is an affine sublattice of $\Z^3$. Note that a one-point subset of
$\Z^t$ is not an affine sublattice of $\Z^t$ since we have excluded
$\{0\}$ from being a lattice.

A subset of $\Z^t$ is called \emph{(Zariski) closed} if it is the
set of common zeros in $\Z^t$ of a system of polynomials in
$t$-variables. For example, the set of integral points on the
$X_1$-axis of the Euclidean plane $\R^2$ is a Zariski closed subset
of $\Z^t$ since the $X_1$-axis is defined by the equation $x_2=0$.
In general, since a hyperplane orthogonal to the $X_i$-axis of
$\R^t$ is defined by an equation of the form $x_i=a_i$ for some
$a_i\in\R$, the set of integral points in such a hyperplane is a
Zariski closed subset of $\Z^t$. It is an easy exercise to show that
a subset of $\Z$ is Zariski closed if and only if it is $\Z$ or
finite.

Let $A$ be a subset of $\Z^t$. A subset $B$ of $A$ is said to be
\emph{(Zariski) dense} in $A$ if every Zariski closed subset of
$\Z^t$ which contains $B$ contains $A$. For example, a subset of
$\Z$ is Zariski dense in $\Z$ if and only if it is infinite.

Following Bourgain-Gamburd-Sarnak \cite{bourgain-gamburd-sarnak}, we
denote by $\mathcal{P}=\{\pm2,\pm3,\pm5,\pm7,\cdots\}$ the set of
integers which are prime numbers up to a sign, and call a point in
$\Z^t\cap\mathcal{P}^t$ a \emph{prime point}. The
Bourgain-Gamburd-Sarnak's formulation of Dickson's conjecture
\cite{bourgain-gamburd-sarnak}, asks when an affine sublattice of
$\Z^t$ has a Zariski dense subset of prime points.

By the prime number theorem, $\mathcal{P}$ is Zariski dense in $\Z$.
Let $q\neq0,b$ be integers. By the prime number theorem in
arithmetic progressions, $(q\Z+b)\cap\mathcal{P}$ is Zariski dense
in $q\Z+b$ if and only if $\gcd(q,b)=1$.

Let $L$ be an affine sublattice of $\Z^t$. If $t>1$ and the
projection of $L$ to some axis of $\Z^t$ is a constant prime point
of $\Z$, then we call $L$ degenerate. For example, the affine
sublattice $\{(-2,3,2n+1)\mid n\in\Z\}$ is a degenerate lattice of
$\Z^3$ since its projection to the first axis is the constant prime
point $-2$ of $\Z$. Let $A$ be the set of numbers $i$ from
$1,2,\cdots,t$ such that the projection of $L$ to the $X_i$-axis of
$\Z^t$ is a constant prime point of $\Z$. It is a proper subset of
$\{1,2,\cdots,t\}$ since $L$ contains more than one point. It is
easy to see that the projection of $L$ to $\prod_{1\leq j\leq
t,j\not\in A}\Z$ is a non-degenerate affine sublattice of
$\prod_{1\leq j\leq t,j\not\in A}\Z$, we call it the non-degenerate
sublattice associated to $L$. For example, the affine sublattice
$2\Z+1$ is a non-degenerate affine sublattice of $\Z$ associated to
the degenerate affine sublattice $\{(-2,3,2n+1)\mid n\in\Z\}$ of
$\Z^3$.

Let $p$ be a prime number, and $\Z_p:=\Z/p\Z$ be the field of $p$
elements. If $L$ is a non-degenerate affine sublattice of $\Z^t$,
then we denote by $L_p$ the image of $L$ under the natural map
$\Z^t\to\Z_p^t$, and call the quantity
$$\alpha_p(L):=
(\frac{p}{p-1})^t\frac{|L_p\cap(\Z_p^{\times})^t|}{|L_p|}$$ the
local factor of $L$ at $p$. If $L$ is a degenerate affine sublattice
of $\Z^t$, and $L^*$ the associated non-degenerate sublattice, then
we call the quantity $\alpha_p(L):=\alpha_p(L^*)$ the local factor
of $L$ at $p$.

\begin{lemma}Let $L$ be an affine sublattice of
$\Z^t$. If $L\cap\mathcal{P}^t$ is Zariski dense in $L$, then all
local factors of $L$ are nonzero.
\end{lemma}

\begin{proof}We may assume that $L$ is a non-degenerate.
Suppose that $L\cap\mathcal{P}^t$ is Zariski dense in $L$, but there
is a prime number $p$ such that $\alpha_p(L)=0$. Then
$L_p\cap(\Z_p^{\times})^t=\emptyset$. So $p|(n_1\cdots n_t)$
whenever $(n_1,\cdots,n_t)\in L$. It follows that if
$(n_1,\cdots,n_t)\in L\cap\mathcal{P}^t$, then $n_i=\pm p$ for some
$i$. That is any $\vec n\in L\cap\mathcal{P}^t$ must lie on one of
the hyperplanes $x_i=\pm p$, $i=1,\cdots,t$. So $L\cap\mathcal{P}^t$
is contained in the union of these hyerplanes, which is Zariski
closed. As $L\cap\mathcal{P}^t$ is Zariski dense in $L$, $L$ is also
contained in that union. It follows that $L$ must be contained in
one of the hyperplanes, contradicting the the non-degenerateness of
$L$. The lemma is proved.
\end{proof}

Justified by that lemma, we call the vanishing of $\alpha_p(L)$ a
local obstruction at $p$ to the Zariski density of
$L\cap\mathcal{P}^t$ in $L$.

Let $q\neq0,b$ be integers. Then the local obstructions to the
Zariski density of prime points in $q\Z+b$ occur precisely at prime
factors of $\gcd(q,b)$. It follows that all the local obstructions
to the Zariski density of prime points in $q\Z+b$ are passed if and
only if $\gcd(q,b)=1$. So, by the prime number theorem in arithmetic
progressions, the set of prime points in $q\Z+b$ is Zariski dense in
$q\Z+b$ if and only if all local obstructions for $q\Z+b$ to have a
Zariski dense subset of prime points are passed.
\begin{conjecture}[Dickson's
conjecture, Bourgain-Gamburd-Sarnak's formulation
\cite{bourgain-gamburd-sarnak}] Let $L$ be an affine sublattice of
$\Z^t$. If all local obstructions to the Zariski density of prime
points in $L$ are passed, then the set of prime points in $L$ is
Zariski dense in $L$.
\end{conjecture}

Let $L$ the sublattice of $\Z^t$ consisting of vectors
$(n_1,\cdots,n_t)$ such that $n_1,\cdots,n_t$ form an arithmetic
progression. Then for every prime number $p$,
$(1,1+p,1+2p,\cdots,1+(t-1)p)$ gives rise to a point in
$L_p\cap(\Z_p^{\times})^t$. It follows that all local obstructions
to the Zariski density of prime points in $L$ are passed. According
to Dickson's conjecture, the prime points in $L$ should be Zariski
dense in $L$. This is now a marvelous theorem of Green-Tao
\cite{green-tao-primeaps}.

\begin{theorem}[Green-Tao
\cite{green-tao-primeaps}]Let $L$ the sublattice of $\Z^t$
consisting of vectors $(n_1,\cdots,n_t)$ such that $n_1,\cdots,n_t$
form an arithmetic progression. Then the set of prime points in $L$
is Zariski dense in $L$.
\end{theorem}

Following Bourgain-Gamburd-Sarnak \cite{bourgain-gamburd-sarnak}, we
define the \emph{von Mangoldt function} $\Lambda$ on $\Z$ by setting
$\Lambda(n)=\log p$ when $\pm n$ is a power of a prime $p$, and
$\Lambda(n)=0$ otherwise. We define the \emph{$t$-dimensional von
Mangoldt function} on $\Z^t$ by setting
$$\Lambda^{\otimes t}(n_1,\ldots,n_t)=\prod_{i=1}^t\Lambda(n_i).$$

Let $L$ be a sublattice of $\Z^t$. The classical Dickson's
conjecture for $L$ predicts that, if $N$ is large,
$K\subseteq[-N,N]^t$ is convex and of volume $\gg N^t$, and
$a\in\Z^t$ is of size $O(N)$ such that all local obstructions to the
Zariski density of prime points in $a+L$ are passed, then the
quantity
$$\frac{1}{|L\cap K|}\sum_{\vec{n}\in L\cap K}\Lambda^{\otimes
t}(a+\vec{n}),$$ which is the expectation of the $t$-dimensional von
Mangoldt function on $a+L\cap K$, is asymptotically
$\prod_p\alpha_p(a+L)$, which, by Lemma 1.3 of
\cite{green-tao-equations}, is always convergent.

\begin{conjecture}[Dickson's conjecture]Let $L$ be a sublattice of $\Z^t$.
If $N$ is large, $K\subseteq[-N,N]^t$ is convex and of volume $\gg
N^t$, and $a\in\Z^t$ is of size $O(N)$ such that all local
obstructions to the Zariski density of prime points in $a+L$ are
passed. Then $$\frac{1}{|L\cap K|}\sum_{\vec{n}\in L\cap
K}\Lambda^{\otimes t}(a+\vec{n})=\prod_p\alpha_p(a+L)+o(1).$$
\end{conjecture}

Green-Tao \cite{green-tao-equations} also made great advances in the
study of Dickson's conjecture for a general sublattice $L$. We
proceed to introduce their main result in this direction.

A map $\Psi$ from $\mathbb{Z}^d$ to $\mathbb{Z}^t$ is called linear
if it is a group homomorphism, is called affine-linear if $\dot
\Psi:=\Psi-\Psi(0)$ is linear. A linear map from $\mathbb{Z}^d$ to
$\mathbb{Z}$ is called a linear form on $\mathbb{Z}^d$. And an
affine-linear map from $\mathbb{Z}^d$ to $\mathbb{Z}$ is called an
affine-linear form on $\mathbb{Z}^d$.

If $\psi_1,\ldots,\psi_t$ is a system of affine-linear forms on
$\mathbb{Z}^d$, then
$$(\psi_1,\ldots,\psi_t)(\vec{n}):=(\psi_1(\vec{n}),\ldots,\psi_t(\vec{n}))$$
is an affine-linear map from $\Z^d$ to $\Z^t$. Conversely, if
$\Psi:\mathbb{Z}^d\to\mathbb{Z}^t$ is an affine-linear map, then
there is a unique system $\psi_1,\ldots,\psi_t$ of affine-linear
forms on $\mathbb{Z}^d$ such that $\Psi=(\psi_1,\ldots,\psi_t)$.
Moreover, $\Psi$ is linear if and only if all $\psi_i$ ($1\leq i\leq
t$) are linear.

If $\Psi:\Z^{d}\to\Z^t$ is a linear map, then the range $\Psi(\Z^d)$
is a sublattice of $\Z^t$. Conversely, if $L$ is a sublattice of
$\Z^t$, then there is a linear map $\Psi:\Z^{d}\to\Z^t$ such that
$L=\Psi(\Z^d)$. Similarly, if $\Psi:\Z^{d}\to\Z^t$ is an
affine-linear map, then the range $\Psi(\Z^d)$ is an affine
sublattice of $\Z^t$. And, if $L$ is an affine sublattice of $\Z^t$,
then there is an affine-linear map $\Psi:\Z^{d}\to\Z^t$ such that
$L=\Psi(\Z^d)$.

\begin{definition}[Green-Tao partition]  Let $\psi_1,\ldots,\psi_t$
be a system of linear forms on $\Z^d$, and $\Psi =
(\psi_1,\ldots,\psi_t)$. A \emph{Green-Tao partition of size
$s\geq0$} of $\Psi$ at $i$ ($1 \leq i \leq t$) is a partition of the
$t-1$ forms $\{ \psi_j: j \neq i\} $ into $s+1$ classes such that
$\psi_i$ does not lie in the $\Q$-linear span of any class.
\end{definition}

\begin{lemma}[Green-Tao, Lemma 1.6 in \cite{green-tao-equations}]  Let $\psi_1,\ldots,\psi_t$
be a system of linear forms on $\Z^d$, and $\Psi =
(\psi_1,\ldots,\psi_t)$. Let $1 \leq i \leq t$. A Green-Tao
partition of $\Psi$ at $i$ exists if and only if $\psi_i$ is not a
rational multiple of any other form.
\end{lemma}

\begin{proof}First suppose that $\psi_i$ is not a
rational multiple of any other form. Then the singleton partition
$\{ \psi_j: j \neq i \}=\cup_{j\neq i}\{\psi_j\}$ is a Green-Tao
partition of $\Psi$ at $i$. Secondly suppose that $\psi_i$ is a
rational multiple of some other form, say $\psi_{i_0}$, and $\{
\psi_j: j \neq i \}=\cup_{k}A_k$ is any partition. Then $\psi_i$
must lie in the $\Q$-linear span of $A_k$ which contains
$\psi_{i_0}$. Therefore that partition is not a Green-Tao partition
of $\Psi$ at $i$.
\end{proof}
\begin{definition}[Green-Tao complexity]  Let $\psi_1,\ldots,\psi_t$
be a system of linear forms on $\Z^d$, and $\Psi =
(\psi_1,\ldots,\psi_t)$. Let $1 \leq i \leq t$. If $\psi_i$ is not a
rational multiple of any other form, the \emph{$i$-complexity} of
$\Psi$ is defined to be the least of the sizes of all Green-Tao
partitions of $\Psi$ at $i$. If $\psi_i$ is a rational multiple of
some other form, the \emph{$i$-complexity} of $\Psi$ is defined to
be $\infty$. The \emph{Green-Tao complexity} of $\Psi$ is defined to
be the largest of all $i$-complexities of $\Psi$ with $1\leq i\leq
t$.
\end{definition}

We reformulate Green-Tao's theorem on Dickson's conjecture in
\cite{green-tao-equations} as follows.
\begin{theorem}[Green-Tao,
\cite{green-tao-equations}]\label{green-tao-theorem-reformulated}
Let $\Psi:\Z^d\to\Z^t$ be a linear map of complexity at most
$s<\infty$. Assume that the Gowers Inverse conjecture $\GI(s)$ and
the \mobname conjecture $\MN(s)$ are true for the parameter $s$.
Then Dickson's conjecture is true for $\Psi(\Z^d)$.
\end{theorem}

The conjectures $\GI(s)$ and $\MN(s)$ were formulated by Green-Tao
\cite{green-tao-equations}. The conjecture $\GI(1)$ is easy.
Green-Tao \cite{green-tao-mobius} reduced the conjecture $\MN(1)$ to
a classical result of Davenport \cite{davenport}. But, according to
Green-Tao \cite{green-tao-equations}, $\MN(1)$ was essentially
already present in the work of Hardy-Littlewood
\cite{hardy-littlewood} and Vinogradov \cite{vinogradov}. The
conjectures $\GI(2)$ and $\MN(2)$ were also proved by Green-Tao
\cite{green-tao-inverse, green-tao-mobius}.

The above reformulation of Green-Tao's theorem on Dickson's
conjecture will be proved in the next section. It needs some
knowledge of the local factors of affine sublattices.

To state the main result of this paper, we introduce the following
definition.

\begin{definition}Let $T:\Z^d\to\Z^d$ be a linear map. We call $T$ \emph{complexity preserving} if for
every linear map $\Psi:\mathbb{Z}^d\to\mathbb{Z}^t$, the Green-Tao
complexities of $\Psi$ and $\Psi\circ T$ are equal.
\end{definition}

The following theorem is the main result of this paper.
\begin{theorem}\label{main-theorem}
Let $\Psi:\Z^d\to\Z^t$ be a linear map, and $L$ be a sublattice of
$\Psi(\Z^d)$ of finite index. Then there is a complexity-preserving
linear map $T:\Z^d\to\Z^d$ such that $T(\Z^d)=\Psi^{-1}(L)$.
\end{theorem}

One can infer the following corollary from Theorems
\ref{main-theorem} and \ref{green-tao-theorem-reformulated}.
\begin{corollary}
Let $\Psi:\Z^d\to\Z^t$ be a linear map of complexity at most
$s<\infty$. Let $L$ be a sublattice of $\Psi(\Z^d)$ of finite index.
Suppose that the Gowers Inverse conjecture $\GI(s)$ and the \mobname
conjecture $\MN(s)$ are true for the parameter $s$. Then Dickson's
conjecture is true for $L$.
\end{corollary}

\begin{proof}By Theorem \ref{main-theorem}, there is a complexity-preserving
linear map $T:\Z^d\to\Z^d$ such that $T(\Z^d)=\Psi^{-1}(L)$. So
$\Psi\circ T$ is of complexity at most $s$. By Theorem
\ref{green-tao-theorem-reformulated}, Dickson's conjecture is true
for $\Psi\circ T(\Z^d)=L$. The corollary is proved.\end{proof}

A great deal of work were done for special sublattices of the range
of the linear map
$$\Psi:\Z^2\to\Z^3,\ (n_1,n_2)\mapsto(n_1,n_2,-n_1-n_2.$$ The pioneer is
Rademacher \cite{rademacher-ternary}. The followers are Ayoub
\cite{ayoub-ternary}, Liu-Tsang \cite{liu-tsang-general}, Liu-Wang
\cite{liu-wang-general}, Liu-Zhan \cite{liu-zhan-ternary}, Bauer
\cite{bauer-ternary}, and etc..

{\it Acknowledgements} The author would like to thank the
Morningside Center of Mathematics, Chinese Academy of Sciences for
support over several years. The author would also like to thank
Jianya Liu for his lecture at the Morningside Center of Mathematics,
which drew the author's attention to some terminologies adopted by
Bourgain-Gamburd-Sarnak \cite{bourgain-gamburd-sarnak}.
\section{Green-Tao's theorem}
Let $p$ be a prime number. Following Green-Tao
\cite{green-tao-equations}, we denote the function
$\frac{p}{p-1}{\bf 1}_{(\Z_p)^{\times}}$ on $\Z_p$ as $\Lambda_p$.
So $\Lambda_p(b)=0$ when $b=p\Z$ and $\Lambda_p(b)=\frac{p}{p-1}$
otherwise. We define the \emph{($t$-dimensional) local von Mangoldt
function} at $p$ on $(\Z_p)^t$ by setting
$$\Lambda_p^{\otimes t}(n_1,\ldots,n_t)=\prod_{i=1}^t\Lambda_p(n_i).$$

In this section we deduce Theorem
\ref{green-tao-theorem-reformulated} from the following theorem,
which is more close to Green-Tao's original formulation.
\begin{theorem}[Green-Tao,
\cite{green-tao-equations}]Let $\Psi:\Z^d\to\Z^t$ be a linear map of
complexity at most $s<\infty$. Assume that the Gowers Inverse
conjecture $\GI(s)$ and the \mobname conjecture $\MN(s)$ are true
for the parameter $s$. If $N$ is large, $K\subseteq[-N,N]^d$ is
convex and of volume $\gg N^d$, and $a\in\Z^t$ is of size $O(N)$
such that all local obstructions to the Zariski density of prime
points in $a+\Psi(\Z^d)$ are passed. Then
$$\frac{1}{|K\cap\Z^d|}\sum_{\vec{n}\in K\cap\Z^d}\Lambda^{\otimes
t}(a+\Psi(\vec n)) =\prod_p\gamma_p(a+\Psi)+o(1),$$ where
$$
\gamma_p(a+\Psi) := \frac{1}{p^d}\sum_{\vec{n} \in \Z_p^d}
\Lambda_p^{\otimes t}(a+\Psi(\vec{n})).$$
\end{theorem}

First we reformulate that theorem as follows.
\begin{theorem}[Green-Tao,
\cite{green-tao-equations}]\label{green-tao-theorem}Let
$\Psi:\Z^d\to\Z^t$ be a linear map of complexity at most $s<\infty$.
Assume that the Gowers Inverse conjecture $\GI(s)$ and the \mobname
conjecture $\MN(s)$ are true for the parameter $s$. If $N$ is large,
$K\subseteq[-N,N]^d$ is convex and of volume $\gg N^d$, and
$a\in\Z^t$ is of size $O(N)$ such that all local obstructions to the
Zariski density of prime points in $a+\Psi(\Z^d)$ are passed. Then
\begin{equation}\label{green-tao-estimate}\frac{1}{|K\cap\Z^d|}\sum_{\vec{n}\in
K\cap\Z^d}\Lambda^{\otimes t}(a+\Psi(\vec n))
=\prod_p\alpha_p(a+L)+o(1),\end{equation} where $L=\Psi(\Z^d)$.
\end{theorem}
\begin{proof}Since $s$ is finite, $a+L$ is non-degenerate.
The theorem then follows from the following one.\end{proof}

\begin{theorem}\label{local-factor}Let
$\Psi:\Z^d\to\Z^t$ be a linear map, and $L$ be a sublattice of
$\Psi(Z^d)$ of finite index. Let $a\in Z^t$ be such that $a+L$ is
non-degenerate. Then
$$
\alpha_p(a+L)=\frac{1}{p^d}\frac{|\Psi(\Z_p^d)|}{|L_p|}\sum_{\vec{n}
\in \Z_p^d,\Psi(\vec n)\in L_p} \Lambda_p^{\otimes
t}(a+\Psi(\vec{n})).$$
\end{theorem}
\begin{proof}Since the
sequence
$$0\to\{\vec{n} \in \Z_p^d \mid
\Psi(\vec
n)=0\}\to\Z_p^d\stackrel{\Psi}{\to}\Z_p^t\to\Z_p^t/\Psi(\Z_p^d)\to0$$
is exact, we have
$$|\{\vec{n} \in \Z_p^d \mid \Psi(\vec n)=0\}|=\frac{p^d}{|\Psi(\Z_p^d)|}.$$
Therefore we have
$$ \sum_{\vec{n} \in \Z_p^d, \Psi(\vec n)\in L_p}
\Lambda_p^{\otimes t}(a+\Psi(\vec{n})) = \sum_{\vec{m} \in L_p}
\Lambda_p^{\otimes t}(\vec{m})\sum_{\vec{n} \in \Z_p^d, \Psi(\vec
n)=\vec m}1$$$$=\sum_{\vec{m} \in L_p} \Lambda_p^{\otimes
t}(\vec{m})\sum_{\vec{n} \in \Z_p^d, \Psi(\vec
n)=0}1=\frac{p^d}{|\Psi(\Z_p^d)|}\sum_{\vec{n} \in L_p}
\Lambda_p^{\otimes t}(\vec{n}).$$ The lemma now follows from the
definition of local factors.
\end{proof}

We now prove Theorem \ref{green-tao-theorem-reformulated}.
\begin{proof}[Proof of Theorem \ref{green-tao-theorem-reformulated}]
Let $r$ be the rank of $\Psi$. Then there is an automorphism $T$ of
$\Z^d$ such that $$\Psi\circ T(n_1,\cdots,n_t)=\Psi\circ
T(n_1,\cdots,n_r,0,\cdots,0), \forall(n_1,\cdots,n_t)\in\Z^d.$$ Let
$K'\subseteq[-N,N]^r$ be convex and of volume $\gg N^r$, and set
$K=K'\times[-N,N]^{d-r}$. Note that \eqref{green-tao-estimate} is
also valid if $\Psi$ is replaced by $\Psi\circ T$. It follows that
\begin{equation}\label{green-tao-estimate}\frac{1}{|K'\cap\Z^r|}\sum_{\vec{n}\in
K'\cap\Z^r}\Lambda^{\otimes t}(a+\Phi(\vec n))
=\prod_p\alpha_p(a+L)+o(1),\end{equation} where $L=\Psi(\Z^d)$, and
$\Phi:\Z^r\to\Z^t$ is defined by setting
$$\Phi(n_1,\cdots,n_r)=\Psi\circ T(n_1,\cdots,n_r,0,\cdots,0).$$ Now
we let $K\subseteq[-N,N]^t$ be convex and of volume $\gg N^t$, and
set $K'=\Phi^{-1}({K})$. Note that there are positive constants
$c_1,c_2$ such that
$$\Phi^{-1}([-N,N]^t)\subseteq[-c_1N,c_2N]^r.$$ It follows that
\eqref{green-tao-estimate} is also valid for this new $K'$. As
$$|K'\cap \Z^d| =|L\cap K|,$$ and
$$\sum_{\vec{n} \in K'\cap \Z^d} \Lambda^{\otimes
t}(a+\Phi(\vec n))=\sum_{\vec{m} \in L\cap K}\Lambda^{\otimes
t}(a+\vec m),$$ \eqref{green-tao-estimate} gives
$$\frac{1}{|L\cap K|}\sum_{\vec{n}\in L\cap K}\Lambda^{\otimes
t}(a+\vec n) =\prod_p\alpha_p(L)+o(1).$$Theorem
\ref{green-tao-theorem-reformulated} is proved.
\end{proof}
\section{Complexity of sublattices of finite index}
In this section we prove Theorem \ref{main-theorem}.
\begin{proof}[Proof of Theorem \ref{main-theorem}]
It is easy to see that $\Psi^{-1}(L)$ is a sublattice of $\Z^d$.
Since $L$ is of finite rank in $\Psi(\Z^d)$, $\Psi^{-1}(L)$ is also
of finite index in $\Z^d$. So $\Psi^{-1}(L)$ can be generated by $d$
vectors in $\Z^d$ which are linearly independent over $\Q$. These
$d$ vectors give rise to a linear map
$T:\mathbb{Z}^d\to\mathbb{Z}^d$ of rank $d$ such that
$T(\Z^d)=\Psi^{-1}(L)$. Theorem \ref{main-theorem} then follows from
the following theorem.
\end{proof}

\begin{theorem}[Complexity Preserving Theorem]\label{complexity-preserving-theorem}Let
$T:\mathbb{Z}^d\to\mathbb{Z}^d$ be a linear map. Then $T$ is
complexity preserving if and only if $T(\Z^d)$ is of rank d.
\end{theorem}

\begin{proof}First we show that $T$ is of rank $d$
under the assumption that $T$ is complexity preserving. Let
$I:\Z^d\to\Z^d$ be the identity map. Since $T$ is complexity
preserving, $T=I\circ T$ is of Green-Tao complexity equal to that of
$I$. It is easy to see that $I$ is of Green-Tao complexity $0$. So
$T$ is of the Green-Tao complexity $0$. Let $T_1,\ldots,T_t$ be
linear forms on $\Z^t$ such that $T=(T_1,\ldots,T_t)$. Then the
system $T_1,\ldots,T_t$ is $\Q$-linearly independent. It follows
that $T(\Z^d)$ is of rank $d$.

We now show that $T$ is complexity preserving under the assumption
that $T$ is of rank $d$. Let $\Psi:\Z^d\to\Z^t$ be an arbitrary
linear map. We must show that the Green-Tao complexity of $\Psi\circ
T$ is equal to that of $\Psi$. By the following lemma, it suffices
to show that the Green-Tao complexity of $\Psi$ is no less than that
of $\Psi\circ T$. As $T(\Z^d)$ is of rank $d$, there is a linear map
$Q:\Z^d\to\Z^d$ and a nonzero integer $a$ such that $T\circ Q=aI$.
Again, by the following lemma, the Green-Tao complexity of
$a\Psi=\Psi\circ T\circ Q$ is no less than that of $\Psi\circ T$.
Since scalars are obviously complexity preserving, the Green-Tao
complexity of $a\Psi$ is equal to that of $\Psi$. It follows that
the Green-Tao complexity of $\Psi$ is no less than that of
$\Psi\circ T$. The proposition is proved.\end{proof}

\begin{lemma}Let $T:\mathbb{Z}^d\to\mathbb{Z}^d$ and
$\Psi:\mathbb{Z}^d\to\mathbb{Z}^t$ be linear. Then the Green-Tao
complexity $\Psi\circ T$ is no less than that of $\Psi$.
\end{lemma}
\begin{proof}Suppose that $\Psi =
(\psi_1,\ldots,\psi_t)$, where $\psi_1,\ldots,\psi_t$ is a system of
linear forms on $\Z^d$. Let $1\leq i\leq t$. Let $s_i$ be the
$i$-complexity of $\Psi\circ T= (\psi_1\circ T,\ldots,\psi_t\circ
T)$. It suffices to show that the $i$-complexity of $\Psi$ is no
greater than $s_i$. We may assume that $s_i<\infty$. Let
$$\{\psi_j\circ T:j\neq i\}=\cup_{k=1}^{s_i}\{\psi_j\circ T:j\in
J_k\}$$ be a Green-Tao partition of $\Psi\circ T$ at $i$. Then
$$\{\psi_j:j\neq i\}=\cup_{k=1}^{s_i}\{\psi_j:j\in J_k\}$$ is a
Green-Tao partition of $\Psi\circ T$ at $i$. It follows that the
$i$-complexity of $\Psi$ is no greater than $s_i$. \end{proof}

\end{document}